\documentclass[leqno]{amsart}
\usepackage{amsmath,amsfonts,amsthm,amssymb, mathrsfs, mathabx, accents, mathdots, textcomp, fontenc,  vntex}
\usepackage[bookmarks]{hyperref}
\hypersetup{backref, colorlinks=true}
\usepackage[french, english]{babel}
\usepackage{amsmath,amsfonts,amsthm,amssymb, mathrsfs, mathabx, accents, mathdots}
\usepackage[all]{xy}
 \usepackage{graphicx}
\usepackage[all]{xy}
\newtheorem{theorem}{Th\'eor\`eme}[section]

\newtheorem{lemma}[theorem]{Lemme}
\newtheorem{proposition}[theorem]{Proposition}

\newtheorem{remark}[theorem]{Remarque}

\newtheorem{example}[theorem]{Exemple}
\newtheorem{definition}[theorem]{D\'efinition}

\newenvironment{equationth}{\stepcounter{theorem}\begin{equation}}{\end{equation}}
\newenvironment{preuve}{{\em{\noindent \textbf{Preuve.} }}}
{\hfill $\blacksquare$}

\def\C{ \mathbb{C}}

\def\R{ \mathbb{R}}

\def\N{ \mathbb{N}}

\thanks{Cet article a \'et\'e \'ecrit avec le support de la R\'egion Provence-Alpes-C\^ote d'Azur (France)  
et  le support de la bourse post-doctorale n$^{{\rm o}}$ 2013/18706-7 de la FAPESP (Br\'esil).}

\def\rond{\mathaccent"7017}

\begin{document}

\large 

\title[]{LA M\'ETHODE DES  FA{\c c}ONS} 
\makeatother

\author[Nguy\~{\^e}n Th\d{i} B\'ich Th\h{u}y]{Nguy\~{\^e}n Th\d{i} B\'ich Th\h{u}y}
\address[{Nguy\~{\^e}n Th\d{i} B\'ich Th\h{u}y}]{UNESP, Universidade Estadual Paulista, ``J\'ulio de Mesquita Filho'', S\~ao Jos\'e do Rio Preto, Brasil}
\email{bichthuy@ibilce.unesp.br}
\maketitle \thispagestyle{empty}
\begin{abstract}
We give a natural method, called {\it ``la m\'ethode des fa{\c c}on''} to stratify the asymptotic variety associated to a dominant polynomial mapping $F: \C^n \to \C^n$. The obtained stratification is differentiable and satisfies the frontier properties. 

\medskip 

\medskip 

\noindent R\'ESUM\'E. Nous donnons une m\'ethode naturelle, que nous appelons {\it ``la m\'ethode des fa{\c c}ons''} pour stratifier la vari\'et\'e  asymptotique associ\'ee \`a une application polynomiale dominante $F: \C^n \to \C^n$. La stratification
obtenue est diff\'erentiable et satisfait la propri\'et\'e de fronti\`ere. 
\end{abstract}

\medskip

\begin{center} 
{\large \bf INTRODUCTION}
\end{center}

\medskip

Soit $F : \C^n_{(x)} \to \C^n_{(\alpha)}$ une application polynomiale. L'ensemble asymptotique de $F$, not\'e $S_F$, est l'ensemble des points du but en lesquels l'application $F$ n'est pas propre. 
Dire que l'application $F$ n'est pas propre en un point  du but $a \in  \C^n_{(\alpha)}$ peut \^etre caract\'eris\'e des deux mani\`eres suivantes~: 

1) Il existe une suite  $\{ \xi_k\}_{k \in \N} \subset \C^n_{(x)}$ dans la source telle que $\{ \xi_k\}$ tende vers l'infini et telle que l'image $F(\xi_k)$ tende vers $a$. 

\noindent Ici, ``la suite $\{ \xi_k\}$ tend vers l'infini'' signifie que la norme euclidienne $\vert \xi_k \vert$ de $\xi_k$ dans $\C^n_{(x)}$ tend vers l'infini. 

2) Il existe une courbe diff\'erentiable  
$$\gamma : (0, +\infty) \, \to \C^n_{(x)}, \quad \gamma(u) = (\gamma_1(u), \ldots, \gamma_n(u))$$ 
tendant vers l'infini et telle que $F \circ \gamma (u)$ tend vers $a$ lorsque $u$ tend vers l'infini.  

Dans cet article, nous travaillons avec la deuxi\`eme caract\'erisation de l'ensemble asymptotique.  C'est-\`a-dire, nous  consid\'erons des courbes diff\'erentiables $\gamma : (0, +\infty) \, \to \C^n_{(x)}$ tendant vers l'infini telles que leurs images $F \circ \gamma (u)$ ne tend pas vers l'infini  lorque $u \in (0, +\infty)$ tend vers l'infini. 
Nous voyons qu'il suffit, pour d\'efinir $S_F$ de consid\'erer des courbes $\gamma$ tendant vers l'infini au sens suivant~: chaque coordonn\'ee $\gamma_1(u), \ldots, \gamma_n(u)$ de cette courbe ou bien tend vers l'infini ou bien converge. C'est ce que nous ferons dans cet article.

Dans les ann\'ees 90, Jelonek a \'etudi\'e l'ensemble asymptotique associ\'e \`a une application polynomiale $F: \C^n \to \C^n$ de mani\`ere approfondie et il en a d\'ecrit les  principales propri\'et\'es \cite{Jelonek1, Jelonek2, Jelonek3}. 
La compr\'ehension de la structure de cet ensemble est tr\`es importante par sa relation avec la  Conjecture  Jacobienne (voir, par exemple, \cite{Essen}).  

Cet article d\'ecrit une m\'ethode naturelle, appel\'ee {\it la m\'ethode des fa{\c c}ons}, pour stra\-ti\-fi\-er la vari\'et\'e asymptotique d'une application polynomiale dominante $F: \C^n \to \C^n$. 
 La d\'etermination des strates de la stratification de l'ensemble asymptotique n\'ecessite plusieurs \'etapes. 
La premi\`ere \'etape de la m\'ethode des fa\c cons fournit une partition dont les \'el\'ements peuvent \^etre des 
vari\'et\'es singuli\`eres. Un raffinement de cette partition est n\'ecessaire, en utilisant une version 
plus fine  de notre m\'ethode (les fa\c cons  ``\'etoile''). 
De mani\`ere plus pr\'ecise, les courbes $\gamma$ tendant vers l'infini et telles que leurs images $F \circ \gamma$ tendent vers les points de $S_F$ seront  labellis\'ees sous la forme de ``{\it fa{\c c}ons}'' et ``{\it fa{\c c}ons \'etoile}'', respectivement, \`a l'aide des deux   \'etapes suivantes~:

A) Dans la premi\`ere \'etape, nous d\'efinissons ``{\it une fa{\c c}on du point $a \in S_F$}''. 
Un point $a$ de $S_F$ est la limite de $ F\circ \gamma (u)$, o\`u $\gamma(u) = (\gamma_1(u), \ldots, \gamma_n(u)) : (0, +\infty) \, \to \C^n_{(x)}$ est une courbe  tendant vers l'infini. Nous classons les 
coordonn\'ees $\gamma_1(u), \ldots, \gamma_n(u)$ de la courbe $\gamma (u)$ en trois cat\'egories~: 
i) les coordonn\'ees $\gamma_{i_r}(u)$ tendant vers l'infini (cette cat\'egorie n'est pas vide);
ii) les coordonn\'ees $\gamma_{j_s}(u)$ telles que $\lim_{u \to \infty} \gamma_{j_s}(u)$ est un nombre complexe  ``ind\'ependant du 
point $a$ dans un voisinage de $a$ dans $S_F$''. 
Cela signifie qu'il existe des points $a'$ voisins de $a$ dans $S_F$ 
et des courbes $\gamma^{a'} = (\gamma^{a'}_1(u), \ldots, \gamma^{a'}_n(u)): (0, +\infty) \, \to \C^n_{(x)}$ tendant vers l'infini telles que $\lim_{u \to \infty} F \circ \gamma^{a'}(u) = a'$ et $\lim_{u \to \infty} \gamma^{a'}_{j_s} (u) = \lim_{u \to \infty} \gamma_{j_s} (u) = constante$; 
iii) les coordonn\'ees $\gamma_{ i_l} (u)$ telles que $\lim_{u \to \infty} \gamma_{i_l}(u)$ est un nombre complexe  ``d\'ependant du point $a$'' (ce cas est le cas contraire du cas ii)).
 L'exemple \ref{exfacon} illustre ces trois  cat\'egories.

Nous d\'efinissons ``{\it une fa{\c c}on du point $a \in S_F$}''
comme un $(p,q)$-uple $(i_1, \ldots , i_p)[j_1, \ldots, j_q]$ d'entiers o\`u  les entiers $i_1, \ldots , i_p$ sont
les indices des coordonn\'ees de la premi\`ere cat\'egorie et $j_1, \ldots, j_q$ 
les indices des coordonn\'ees de la seconde cat\'egorie (D\'efinition \ref{definitionXi}). 
 La relation d'\'equivalence $a_1  \sim a_2$ si et seulement si les ensembles de fa{\c c}ons de $a_1$ et $a_2$ 
 co\"\i ncident d\'etermine une partition finie de $S_F$. 
 La partition obtenue dans cette premi\`ere \'etape n'est en fait pas compos\'ee de parties lisses. 
 
 B) Dans la seconde \'etape,  nous d\'eterminons la  partie singuli\`ere et  la  partie r\'eguli\`ere de chaque \'el\'ement   d\'ecrite dans l'\'etape A.  Nous distinguons les comportements des courbes correspondant aux  diff\'erentes fa{\c c}ons. Ceux-ci  seront formalis\'ees sous la forme ``{fa{\c c}ons \'etoile}''~: chaque fa{\c c}on ``\'etoile'' diff\'erente d\'efinit des  courbes co\-rres\-pon\-dantes parall\`eles, qui, en fait, d\'efinissent un feuilletage ({\it cf.} Exemple \ref{exfeuilletage}). Les images des feuilletages  co\-rres\-pon\-dant aux  diff\'erentes fa{\c c}ons d'approcher la partie  singuli\`ere de $S_F$, ce qui nous permet de d\'ecomposer celle-ci en strates. 
 
Nous obtenons une stratification diff\'erentiable de l'ensemble asymptotique $S_F$ (Th\'eor\`eme \ref{theostraXi*}). Cette stratification satisfait aussi la propri\'et\'e de fronti\`ere.

Une autre motivation de cette \'etude r\'eside dans le fait que conna\^itre une stratification d'une vari\'et\'e singuli\`ere permet de calculer son homologie d'intersection~: Dans \cite{Valette}, les  auteurs ont contruit une vari\'et\'e singuli\`ere $V_F$ associ\'ee \`a une application polynomiale  $F: \C^2 \to \C^2$ telle que son homologie d'intersection caract\'erise  la propret\'e de $F$ dans le cas o\`u le jacobien de $F$ est partout non nul. L'article \cite{Valette} fournit donc une nouvelle approximation de la Conjecture Jacobienne via l'homologie d'intersection de la vari\'et\'e singuli\`ere $V_F$. 
De plus, une stratification de l'ensemble asymptotique $S_F$ permet de fournir une stratification de la vari\'et\'e $V_F$. Nous pr\'ecisons donc les r\'esultats obtenus dans \cite{Valette} et aussi dans \cite{Thuy2} (\cite{Thuy2} est une g\'en\'elisation de \cite{Valette}).

Comme {\it la m\'ethode des fa{\c c}ons} est une m\'ethode nouvelle, cet article est agr\'ement\'e  de nombreux  exemples afin de mieux comprendre l'id\'ee de la m\'ethode, ainsi que la structure de l'ensemble asymptotique associ\'e \`a une application polynomiale $F: \C^n \to \C^n$. 

\subsection*{Remerciements}  C'est un grand plaisir de remercier J.-P. Brasselet pour son int\'er\^et et ses encouragements.

\medskip 

\begin{center}
{ \large \bf NOTATION}
\end{center}

\medskip


Nous consid\'erons dans cet article des applications polynomiales $F: \C^n \to \C^n$. Nous \'ecrivons souvent $F: \C^n_{(x)} \to \C^n_{(\alpha)}$ pour distinguer entre la source et le but. 

\medskip

\medskip 

\begin{center}
{ \large \bf 0. PR\'ELIMINAIRES}
\end{center}

\medskip

\subsection{Stratification} \label{sectionStratification}

\begin{definition}
{\rm 
 {Soit $V$ une vari\'et\'e (diff\'erentiable ou alg\'ebrique, ou analytique) de dimension $m$. {\it Une stratification}  (${{\mathscr{S}}}$) de $V$ est la donn\'ee d'une filtration 
$$V=V_m \supseteq V_{m-1} \supseteq V_{m-2} \supseteq \dots \supseteq V_1 \supseteq V_0 \supseteq V_{-1} = \emptyset $$
de $V$ telle que toutes les diff\'erences $X_{i} = V_{i} \setminus V_{i-1}$ sont ou bien   vides ou bien unions localement finies de sous-vari\'et\'es  lisses connexes et localement ferm\'ees de dimension $i$, appel\'ees strates.


Soit $S_i$ une strate de $V$ et soit $\overline{S_i}$ son adh\'erence dans $V$. Pour toute strate $S_i$ de $V$, 
si $\overline{S_i} \setminus S_i$ est l'union de strates de $V$, 
alors nous disons que la stratification de $V$ satisfait la propri\'et\'e de fronti\`ere. 
}
}
\end{definition}

\subsection{L'ensemble asymptotique.} \label{ensembleJelonek}

Soit $F : \C^n_{(x)} \to \C^n_{(\alpha)}$ une application polynomiale. Notons $S_F$ l'ensemble des points du but pour lesquels l'application $F$ n'est pas propre, {\it i.e.},  
$$S_F = \{ a \in \C^n_{(\alpha)} \text{ tel que } \exists \{ \xi_k\}_{k \in \N} \subset \C^n_{(x)}, \vert \xi_k \vert  \text{ tend vers l'infini et } F(\xi_k) \text{ tend vers } a\},$$
o\`u $ \vert \xi_k \vert$ est  la norme euclidienne de $\xi_k$ dans $\C^n$. 
L'ensemble $S_F$ est appel\'e l'ensemble asymptotique de $F$.

\begin{lemma} \label{lemmecourbe} 
{\rm Le point $a$ appartient \`a $S_F$ si et seulement s'il existe une courbe diff\'erentiable  $\gamma(u) : (0, +\infty) \, \to \C^n_{(x)}$ tendant vers l'infini et telle que $F \circ \gamma (u)$ tend vers $a$ lorsque $u$ tend vers l'infini. 
}
\end{lemma}

Rappelons qu'il suffit, pour d\'efinir $S_F$ de consid\'erer des courbes $\gamma(u) = (\gamma_1(u), \ldots, \gamma_n(u))$ tendant vers l'infini, au sens suivant ~: chaque coordonn\'ee $\gamma_1(u), \ldots, \gamma_n(u)$ de cette courbe ou bien tend vers l'infini ou bien converge. C'est ce que nous ferons dans cet article.

\begin{definition} \label{definitiondominant}
{\rm Une application polynomiale $F : \C^n_{(x)} \to \C^n_{(\alpha)}$ est dite {\it dominante} si l'adh\'erence de $F(\C^n_{(x)})$ est dense dans $\C^n_{(\alpha)}$, 
c'est-\`a-dire $\overline{F(\C^n_{(x)})} = \C^n_{(\alpha)}$. 
}
\end{definition}

L'ensemble $S_F$ a \'et\'e intens\'ement \'etudi\'e par Jelonek dans une s\'erie d'articles \cite{Jelonek1, Jelonek2, Jelonek3}. En particulier, il a montr\'e le r\'esultat suivant~:

\begin{theorem} \cite{Jelonek1} \label{theoremjelonek1}
Soit $F = (F_1, \ldots, F_n) : \C^n \rightarrow \C^n$ une application polynomiale dominante. Alors, l'ensemble asymptotique $S_F$ de $F$ ou bien vide, ou bien une hypersurface de l'espace but. 
\end{theorem}

\medskip 

\medskip 

Soit $F: \C^n_{(x)} \to \C^n_{(\alpha)}$ une application polynomiale dominante. 
{\it La m\'ethode des fa{\c c}ons} consiste \`a classer les coordonn\'ees des courbes  diff\'erentiables  $\gamma(u) = (\gamma_1(u), \ldots, \gamma_n(u)) : (0, +\infty) \, \to \C^n_{(x)}$ tendant vers l'infini 
et telles que $F \circ \gamma (u)$ tend vers un point $a$ de l'ensemble asymptotique $S_F$ en trois cat\'egories (D\'efinition \ref{definitionXi}). Nous d\'efinissons alors, 
dans une premi\`ere \'etape, une partition de $S_F$ en parties qui, bien que tr\`es significatives  pour notre construction, se r\'ev\`elent malheureusement  
pouvoir \^etre singuli\`eres (Proposition \ref{profaconetoile}). Un raffinement de la d\'efinition des  fa{\c c}ons (que nous appelons
 fa{\c c}ons ``\'etoile'') est n\'ecessaire afin d'obtenir une stratification en strates lisses 
 (D\'efinition  \ref{defpreetoile} et Th\'eor\`eme \ref{theostraXi*}). Prenons d'abord un exemple afin de justifier la premi\`ere \'etape de la m\'ethode qui suit.

\section{Construction des fa{\c c}ons}

\subsection{Un exemple}
\begin{example} \label{exfacon}
{\rm

Soit $F =  (F_1, F_2, F_3) : \C^3_{(x_1, x_2, x_3)} \to \C^3_{(\alpha_1, \alpha_2, \alpha_3)}$ l'application polynomiale dominante telle que 
$$F_1:=x_1, \qquad F_2:= x_2, \qquad F_3:=x_1x_2x_3.$$

Nous d\'eterminons maintenant l'ensemble asymptotique $S_F$ utilisant le Lemme \ref{lemmecourbe}. 
Sup\-po\-sons qu'il existe une courbe $\gamma = (\gamma_1, \gamma_2, \gamma_3) : (0, +\infty) \, \to \C^3_{(x_1, x_2, x_3)}$ tendant vers l'infini telle que $F \circ \gamma (u)$ ne tende pas vers l'infini. 
Comme $F \circ \gamma (u) = (\gamma_1(u), \gamma_2(u), \gamma_1(u) \gamma_2(u) \gamma_3(u))$, 
donc $\gamma_1(u)$ et $\gamma_2 (u)$ ne peuvent pas tendre vers l'infini. 
Alors, nous avons les trois cas sui\-vants :

i) $\gamma_1 (u)$ tend vers 0, $\gamma_2(u)$ tend vers un nombre complexe $\alpha_2 \in \C$ et $\gamma_3 (u)$ tend vers l'infini.   Afin de d\'eterminer le plus grand sous-ensemble possible 
de $S_F$, nous chosissons $\gamma_1 (u)$ tend vers 0 et $\gamma_3 (u)$ tend vers l'infini de telle mani\`ere que le produit $\gamma_1(u) \gamma_3(u)$ tende vers un nombre complexe non nul. 
 Nous choisissons, par exemple  
$$ \gamma(u) = \left( \frac{1}{u}, \alpha_2, \frac{\alpha_3}{\alpha_2}u \right)$$ 
o\`u $\alpha_2 \neq 0$, alors $F \circ \gamma (u)$ tend vers le point $a = (0, \alpha_2, \alpha_3)$ dans $S_F$. 
Nous obtenons la composante $S_2 = \{\alpha_1 = 0\} \setminus 0\alpha_3$, de dimension 2 de $S_F$. 
Nous disons qu'une  {\it ``fa{\c c}on''} de $S_2$ est  $(3)[1]$, o\`u
\begin{enumerate}
\item Le symbole ``(3)''  signifie que la troisi\`eme coordonn\'ee $\gamma_3(u) = \frac{\alpha_3}{\alpha_2}u$ de la courbe $\gamma(u)$ tend vers l'infini. 
\item Le symbole ``[1]'' signifie que  
la premi\`ere coordonn\'ee $\gamma_1(u) = \frac{1}{u}$ de la courbe $\gamma(u)$ tend vers z\'ero, qui est un nombre complexe fix\'e, lequel ne d\'epend pas du point $a = (0, \alpha_2, \alpha_3)$ quand $a$ d\'ecrit  
$S_2 = \{\alpha_1 = 0\} \setminus 0\alpha_3$. 
\item La deuxi\`eme coordonn\'ee $\gamma_2(u) =  \alpha_2$ de la courbe $\gamma(u)$ tend vers un nombre complexe 
 $\alpha_2$ d\'ependant du point  $a = (0, \alpha_2, \alpha_3)$ quand $a$ varie. 
 Alors, avec notre choix de labeliser les courbes, l'indice ``2'' n'apparaitra pas dans la façon (3)[1]
\end{enumerate}

Notons que  nous pouvons v\'erifier facilement que toute courbe $\hat{\gamma} (u)$ tendant vers l'infini telle que $F \circ \hat{\gamma} (u)$ tende vers un point de $S_2$ admet la m\^eme fa{\c c}on $(3)[1]$. Nous disons que  la fa{\c c}on de $S_2 $ est $(3)[1]$. 


ii) $\gamma_1 (u)$ tend vers un nombre complexe $\alpha_1 \in \C$, $\gamma_2(u)$ tend vers z\'ero et  $\gamma_3 (u)$ tend vers l'infini : nous choisissons, par exemple, la courbe 
$$\gamma(u) = \left(\alpha_1, \frac{1}{u}, \frac{\alpha_3}{\alpha_1}u \right),$$
 alors $F \circ \gamma (u)$ tend vers le point $(\alpha_1, 0, \alpha_3)$ de $S_F$. 
De la même manière que dans le cas i), la {\it fa{\c c}on} $(3)[2]$ d\'etermine donc la composante  $S'_2 = \{\alpha_2 = 0\} \setminus 0\alpha_3$ de dimension 2 de $S_F$. 

iii) $\gamma_1 (u)$ et  $\gamma_2(u)$ tendent vers z\'ero, et  $\gamma_3 (u)$ tend vers l'infini:  nous choisissons, par exemple, la courbe
$$ \gamma(u) = \left( \frac{1}{u}, \frac{1}{u}, \alpha_3 u^2 \right),$$
alors $F \circ \gamma (u)$ tend vers le point $(0, 0, \alpha_3)$ de $S_F$. 
 De la même manière que dans le cas i), la  {\it fa{\c c}on} $(3)[1, 2]$ d\'etermine donc la  composante $S_1 = 0\alpha_3$ de dimension 1 de $S_F$.

\medskip

Dans cet exemple illustrant l'id\'ee de  la m\'ethode des fa{\c c}ons, nous 
subdivisons  l'ensemble asymptotique de $F$ en strates lisses. Les strates sont 
$S_2 = \{\alpha_1 = 0\} \setminus 0\alpha_3$, $S'_2 = \{\alpha_2 = 0\} \setminus 0\alpha_3$ et $S_1 = 0\alpha_3$ d\'efinies par les fa{\c c}ons (3)[1], (3)[2] et (3)[1,2], respectivement.
}
\end{example}

 Revenons au cas g\'en\'eral et explicitons la m\'ethode des fa\c cons.

\subsection{Partition de $S_F$ d\'efinie par les  fa{\c c}ons}


\begin{definition} \label{definitionXi}
{\rm 
Soit $F : \C^n_{(x)} \to \C^n_{(\alpha)}$ une application polynomiale dominante telle que 
$S_F \neq \emptyset$.  
Par le lemme \ref{lemmecourbe}, pour chaque point $a$ de $S_F$, il existe une courbe 
$$\gamma^a(u) = (\gamma_1^a(u), \ldots, \gamma_n^a(u)) :  (0, +\infty) \, \to \C^n_{(x)}$$
tendant vers l'infini telle que $F \circ \gamma^a(u)$ tende vers $a$ lorsque $u$ tend vers l'infini. 
 Pour cette courbe, il existe au moins un indice 
$i \in \N \setminus \{0\}$ tel que $\gamma_i^a(u)$ tende  vers l'infini quand $u$ tend vers l'infini. D\'efinissons {\it la fa{\c c}on de tendre vers l'infini de la courbe $\gamma^a$} comme un $(p,q)$-uple maximal $\kappa = (i_1, \ldots , i_p)[j_1, \ldots, j_q]$ d'entiers tous diff\'erents et compris entre 1 et $n$,  tels que 
\begin{enumerate}
\item[i)] $\gamma_{i_r}^{a}(u)$  tend vers l'infini pour tout  $r = 1, \ldots , p$, 

\item[ii)] pour tous $s = 1, \ldots , q$, la courbe $\gamma_{j_s}^a(u)$  tend vers une constante (complexe)  ind\'ependante du point $a$ quand $a$ varie. Cela signifie que~:
\begin{enumerate}
\item[ii.1)] Ou bien il existe dans $S_F$ une sous-vari\'et\'e $U_a$ contenant  $a$ 
telle que pour tout point $a'$ de $U_a$, 
il existe une courbe 
$$\gamma^{a'}(u) = (\gamma_1^{a'}(u), \ldots, \gamma_n^{a'}(u)) :  (0, +\infty) \, \to \C^n_{(x)}$$
tendant vers l'infini telle que 
\begin{enumerate}
\item[a)] $F(\gamma^{a'}(u))$ tend vers $a'$,
\item[b)] $\gamma_{i_r}^{a'}(u)$  tend vers l'infini pour tout  $r = 1, \ldots , p$, 
\item[c)] pour tous $s = 1, \ldots , q$,  
$\lim_{u \to \infty} \gamma_{j_s}^{a'}(u) = \lim_{u \to \infty} \gamma_{j_s}^{a}(u)$ et cette limite est finie.
\end{enumerate}

\item[ii.2)] Ou bien il n'existe pas de telle sous-vari\'et\'e $U_a$ comme dans ii.1) et dans ce cas, nous d\'efinissons 
$$\kappa = (i_1, \ldots , i_p)[j_1, \ldots, j_{n-p}],$$ 
o\`u $\gamma_{i_r}^{a}$  tend vers l'infini pour tout  $r = 1, \ldots , p$ et $\{i_1, \ldots , i_p\} \cup \{j_1, \ldots, j_{n-p}\} = \{1, \ldots, n\}$. Dans ce cas, l'ensemble des points $a$ est  une sous-vari\'et\'e de dimension 0 de $S_F$.
\end{enumerate}
\end{enumerate}

L'application $F: \C^3 \to \C^3$ d\'efinie par $F(x_1, x_2, x_3) = (x_1x_2, x_2x_3, x_3x_1)$ fournit un exemple d'application pour laquelle il existe dans $S_F$ une telle sous-vari\'et\'e de dimension 0 d\'efinie par les trois fa{\c c}ons (1)[2, 3], (2)[1, 3] et (3)[1, 2].

La  fa{\c c}on de tendre vers l'infini de la courbe $\gamma^a$ sera appel\'ee aussi {\it une  
  fa{\c c}on du point $a$} ou {\it une fa{\c c}on de $S_F$}. 

} 
\end{definition}

 Notons qu'en g\'en\'eral, l'ensemble des indices $\{i_1, \ldots , i_p\}$ n'est jamais vide, mais l'ensemble des indices $\{j_1, \ldots , j_q\}$ peut \^etre vide. Notons aussi que 
$$\{i_1, \ldots , i_p\} \cup \{j_1, \ldots , j_q\} \subset \{1, \ldots, n \},$$
 et nous pouvons avoir $\{i_1, \ldots , i_p\} \cup \{j_1, \ldots , j_q\} \neq \{1, \ldots, n \}$ (voir Exemple \ref{exfacon}). 

Dans la suite, nous montrons que la relation d'\'equivalence $a_1  \sim a_2$ si et seulement si les ensembles de fa{\c c}ons de $a_1$ et $a_2$ 
 co\"\i ncident d\'etermine une partition finie de $S_F$. Pour cela, nous devons montrer d'abord que le nombre des fa{\c c}ons possibles de $S_F$ est fini.

\begin{proposition} \label{pro nombre}
{\rm Soit $F : \C^n \to \C^n$ une application polynomiale dominante telle que $S_F \neq \emptyset$. Alors le nombre des fa{\c c}ons possibles de $S_F$ est fini. De mani\`ere pr\'ecise, le nombre maximum des fa{\c c}ons de $S_F$ est \'egal \`a 
$$\sum_{k = 1}^n C_k^n + \sum_{k = 1}^{n-1} C_k^n + \sum_{k = 2}^{n-1} A_k^{n},$$
o\`u
$$C_k^n = \frac{n!}{k!(n-k)!},  \quad \quad \quad A_k^n = \frac{n!}{(n-k)!}.$$
}
\end{proposition}
\begin{preuve}
Soit $\kappa = (i_1, \ldots , i_p)[j_1, \ldots, j_q]$  une fa{\c c}on de l'ensemble asymptotique $S_F$. D'apr\`es la D\'efinition \ref{definitionXi}, nous avons les cas possibles suivants~:

\medskip 

i) Si  $\{ i_1, \ldots, i_p\} \cup \{ j_1, \ldots, j_q\} = \{ 1, \ldots , n\}$ : Nous avons $\sum_{k = 1}^n C_k^n$ fa{\c c}ons  possibles.

\medskip 

ii) Si  $\{ i_1, \ldots, i_p\} \cup \{ j_1, \ldots, j_q\} \neq \{ 1, \ldots , n\}$ et $\{ j_1, \ldots, j_q\} = \emptyset$~: Nous avons $\sum_{k = 1}^{n-1} C_k^n$ fa{\c c}ons  possibles. 

\medskip 

iii) Si  $\{ i_1, \ldots, i_p\} \cup \{ j_1, \ldots, j_q\} \neq \{ 1, \ldots , n\}$ et $\{ j_1, \ldots, j_q\} \neq \emptyset$~: Nous avons $\sum_{k = 2}^{n-1} A_k^{n}$ fa{\c c}ons possibles. 

\medskip  

Donc le nombre maximum des fa{\c c}ons possibles de $S_F$ est  
$$\sum_{k = 1}^n C_k^n + \sum_{k = 1}^{n-1} C_k^n + \sum_{k = 2}^{n-1} A_k^n.$$
\end{preuve}

\begin{example}
{\rm Le nombre maximum des fa{\c c}ons d'une application polynomiale $F: \C^3 \to \C^3$ est 19.}
\end{example}

\begin{proposition} \label{pro partitionfini1}
{\rm Notons $\Xi(a)$ l'ensemble de toutes les fa{\c c}ons du point $a$. La partition de $S_F$ d\'efinie par la relation 
\begin{equationth} \label{partition facon}
a_1 \sim a_2 \text{ si et seulement si } \Xi(a_1) = \Xi(a_2)
\end{equationth}
est une partition finie de $S_F$. 
}
\end{proposition}

\begin{preuve}
Le r\'esultat d\'ecoule de la Proposition \ref{pro nombre}.
\end{preuve}

\begin{definition} \label{Deffacon}
{\rm  La partition de $S_F$ d\'efinie par la relation (\ref{partition facon}) est appel\'ee   {\it partition de $S_F$ d\'efinie par les fa{\c c}ons}. 
}
\end{definition}

\begin{example}
{\rm Revenons \`a l'exemple \ref{exfacon}, la partition de $S_F$ d\'efinie par les fa{\c c}ons  comporte trois sous-vari\'et\'es lisses~:

i) $S_2 = \{\alpha_1 = 0\} \setminus 0\alpha_3$ de dimension 2 et d\'efinie par la fa{\c c}on (3)[1], 

ii) $S'_2 = \{\alpha_2 = 0\} \setminus 0\alpha_3$ de dimension 2 et d\'efinie par la fa{\c c}on (3)[2], 

iii) $S_1 = 0\alpha_3$ de dimension 1 et  d\'efinie par la fa{\c c}on (3)[1,2].
}
\end{example}

Notons qu'en un point d'une sous-vari\'et\'e de dimension maximun de la partition de $S_F$ d\'efinie par les fa{\c c}ons, il est possible d'avoir plusieurs fa{\c c}ons comme le montre l'exemple suivant~: 
$$F: \C^3_{(x_1, x_2, x_3)} \to \C^3_{(\alpha_1, \alpha_2, \alpha_3)}, \quad F(x_1, x_2, x_3) = (x_1x_2 + x_3, x_2x_3 + x_3x_1, x_1x_2 + x_2x_3 + x_3x_1)$$
(l'ensemble $S_F$ n'admet qu'une sous-vari\'et\'e lisse : le plan $\alpha_3 = \alpha_1 + \alpha_2$ qui consiste des deux fa{\c c}ons $(1)[2, 3]$ et $(2)[1, 3]$).

Les sous-vari\'et\'es de la partition de l'ensemble asymptotique $S_F$ d\'efinie par les fa{\c c}ons,  se r\'ev\`elent  malheureusement 
pouvoir \^etre singuli\`eres~: 

\begin{proposition} \label{profaconetoile}
{\rm La partition de l'ensemble asymptotique $S_F$ d\'efinie par les fa{\c c}ons 
 n'est pas une stratification diff\'erentiable. 
}
\end{proposition}

\begin{preuve}
 Soit $F : \C^2_{(x_1, x_2)} \to \C^2_{(\alpha_1, \alpha_2)}$ l'application polynomiale dominante telle que 
$$F(x_1,x_2) = \left({(x_1x_2)}^2, {(x_1x_2)}^3 + x_1 \right).$$ 
Si la courbe $\gamma(u) = (\gamma_1(u), \gamma_2(u)) :  (0, +\infty) \rightarrow  \C^2_{(x_1, x_2)} $ 
tend vers l'infini telle que 
$F \circ \gamma (u) = \left( { (\gamma_{1}(u) \gamma_2 (u))}^2, {(\gamma_{1}(u) \gamma_{2}(u))}^3 + \gamma_{1}(u) \right)$ 
ne tend pas vers l'infini,  alors $\gamma_{1}(u)$ ne peut pas tendre vers l'infini. 
Comme $\gamma(u)$ tend vers l'infini, alors $\gamma_ 2 (u)$ doit tendre vers l'infini et donc $S_F$ n'admet  qu'une seule fa{\c c}on $\kappa = (2)[1]$. Si nous choisissons les courbes coordonn\'ees $\gamma_1(u)$ et $\gamma_2(u)$ tendant vers 0 et l'infini, respectivement, telles que le produit $\gamma_1(u) \gamma_2(u)$ tend vers un nombre complexe $\alpha \in \C$, alors $F\circ \gamma (u)$ tend vers  $(\alpha^2, \alpha^3)$. 
 L'ensemble $S_F$  est donc la courbe $ \alpha_1^3 = \alpha_2^2$ dans $\C^2_{(\alpha_1, \alpha_2)}$ ayant un point singulier \`a l'origine (voir Figure \ref{facon1}). 

\begin{figure}[h!]
\centering
\includegraphics[scale=0.7]{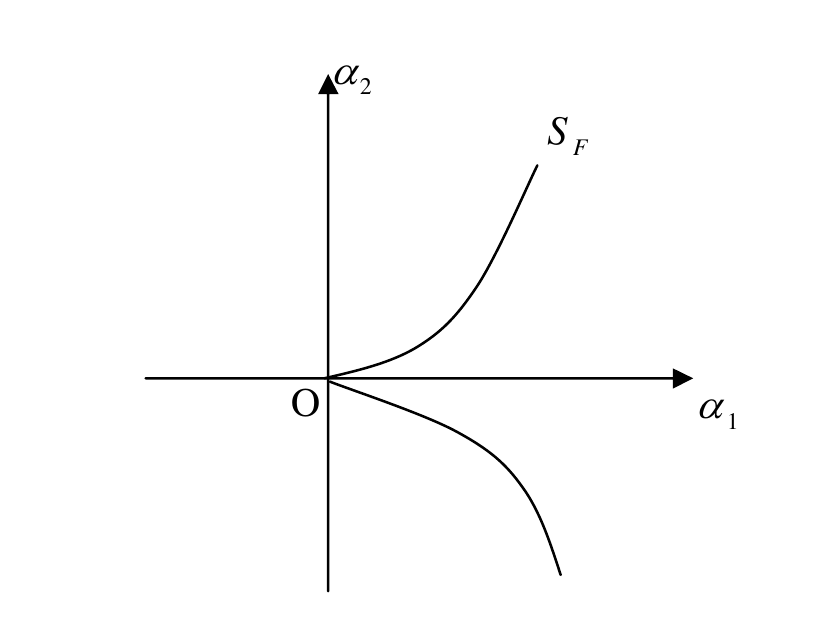}
\caption{L'ensemble asymptotique $S_F$ de l'application $F(x_1,x_2) = \left({(x_1x_2)}^2, {(x_1x_2)}^3 + x_1 \right).$}
\label{facon1}
\end{figure}

\end{preuve}

Donc, un raffinement de la d\'efinition des  fa{\c c}ons  est n\'ecessaire afin d'obtenir une stra\-ti\-fi\-ca\-tion en strates lisses. C'est ce que nous faisons dans la section suivante, en d\'efinissant les {\it fa{\c c}ons \'etoile}.

\section{Construction des  fa{\c c}ons \'etoile}

La question est donc de savoir pourquoi, dans l'exemple de la Proposition \ref{profaconetoile}
$$F : \C^2_{(x_1, x_2)} \to \C^2_{(\alpha_1, \alpha_2)}, 
\quad F(x_1,x_2) = \left({(x_1x_2)}^2, {(x_1x_2)}^3 + x_1 \right),$$ 
la partition de $S_F$ d\'efinie par les fa{\c c}ons n'est pas  une stratification. Nous devons comprendre la diff\'erence entre ``la fa{\c c}on $[2](1)$ \`a l'origine'' et ``la fa{\c c}on $[2](1)$ au point $a \in S_F \setminus \{ 0 \}$''~:
Pour un point $a = (\alpha^2, \alpha^3)$ dans $S_F \setminus \{ 0 \}$, nous devons choisir la  courbe $ \left( \frac{1}{u^r}, \alpha u^s \right)$, 
o\`u $r$ doit \^etre \'egal \`a $s$.  \`A l'origine 0, nous devons choisir une courbe $ \left( \frac{1}{u^r},  u^s \right)$ o\`u $r > s$.  Nous voyons que nous devons \'etudier la vitesse de  tendre vers zero et  de tendre vers l'infini des courbes coordonn\'ees. Dans l'exemple ci-dessus, si la vitesse de tendre vers zero  de la premi\`ere courbe coordonn\'ee est plus rapide que la vitesse de tendre vers l'infini  de la deuxi\`eme courbe coordonn\'ee, nous tendons \`a la partie singuli\`ere de $S_F$.

Afin de formaliser l'id\'ee des {\it fa{\c c}ons \'etoile}, nous avons donc besoin d'abord de d\'ecrire  ``la vitesse de  tendre vers zero et  de tendre vers l'infini des courbes coordonn\'ees dans $\C$''.

\subsection{Quelques d\'efinitions nec\'essaires} 

Notons que si une courbe $\rho :  (0, +\infty) \, \to \C$
 tend vers un nombre complexe $\lambda$ dans $\C$, alors 
nous pouvons consid\'erer que cette courbe tend vers $0$ par changement de  variables $\rho - \lambda$. Pour cette raison, nous pouvons supposer qu'une courbe $\rho$ dans $\C$ ou bien tend vers zero, ou bien tend vers l'infini. 
 Avant de  formaliser l'id\'ee de ``fa{\c c}on \'etoile'', nous donnons les d\'efinitions suivantes~:

\begin{definition} 
{\rm Consid\'erons une courbe $\rho(u):  (0, +\infty) \, \to \C$.  
\begin{enumerate}
\item  Si $\rho(u) = {(\lambda/u)}^t + \cdots$, o\`u $\lambda \in \C \setminus \{ 0 \}$, $t \in \R, t > 0$ et les \'el\'ements dans ``$\ldots$'' sont de la forme ${(\lambda'/u)}^r$ avec $r > t $, 
nous disons que la courbe $\rho$ tend vers 0 avec le degr\'e $t$. 
\item Si $\rho (u) =\lambda u^t + \cdots$, 
o\`u  $\lambda \in \C \setminus \{ 0 \}$, 
$t \in \R, t > 0$ et les \'el\'ements dans ``$\ldots$'' sont de la forme 
$\lambda' k^r$ avec $ r < t$, 
 nous disons que la courbe $\rho$ tend vers l'infini avec le degr\'e $t.$
\end{enumerate}

}
\end{definition}

\begin{definition} \label{uple}
{\rm  Soit   $\kappa = (i_1, \ldots, i_p)[j_1, \ldots, j_q]$ une fa{\c c}on de $S_F$ ({\it cf.} D\'efinition \ref{definitionXi}) et soit $\gamma = (\gamma_1, \ldots, \gamma_n) :  (0, +\infty) \, \to \C^n_{(x)}$  une courbe tendant vers l'infini avec la 
fa{\c c}on $\kappa$.  Supposons que~: 

$\quad \quad \quad$ $\gamma_{i_r}(u)$ tend vers l'infini avec le degr\'e $l_{i_r}$, pour $r=1, \ldots, p$,

$\quad \quad \quad$ $\gamma_{j_s}(u)$ tend vers 0 avec le degr\'e $l_{j_s}$, pour $s=1, \ldots, q$. 

\noindent Le $(p+q)$-uple $(l_{i_1}, \ldots, l_{i_p}, l_{j_1}, \ldots, l_{j_q}) $ est appel\'e {\it $(p+q)$-uple associ\'e \`a la courbe}  
$\gamma$. 
}
\end{definition}

\begin{remark}
{\rm  Par un changement de param\`etre, nous povons toujours supposer que les nombres $l_{i_1}, \ldots, l_{i_p}, l_{j_1}, \ldots, l_{j_q} $ du  $(p+q)$-uple  d'une courbe $\gamma$ dans la d\'efinition \ref{uple} sont des entiers naturels.
}  

\end{remark}

La d\'efinition importante suivante  aide \`a formaliser les {\it fa{\c c}ons \'etoile} par la suite.

\begin{definition} \label{defs:pquple}
{\rm Soit $S_\nu$ une sous-vari\'et\'e de $S_F$ et $\kappa = (i_1, \ldots, i_p)[j_1, \ldots, j_q]$ une fa{\c c}on de $S_\nu$.
Supposons    $\gamma, \gamma': (0, +\infty) \, \to \C^n_{(x)}$ deux courbes tendant vers l'infini avec la m\^eme fa\c con 
$\kappa$  telles que $F \circ \gamma$ et $F \circ \gamma'$ tendent vers deux points de $S_\nu$. 
 Notons  $(l_{i_1}, \ldots, l_{i_p}, l_{j_1}, \ldots, l_{j_q})$ et $(l'_{i_1}, \ldots, l'_{i_p}, l'_{j_1}, \ldots, l'_{j_q})$ leurs  deux $(p+q)$-uples associ\'es aux courbes $\gamma$ et $\gamma'$, respectivement. 
Nous disons que les deux courbes $\gamma$ et $\gamma'$ sont { \it  \'equivalentes}  si leur deux   $(p+q)$-uples associ\'es sont proportionnels, c'est-\`a-dire
\begin{equationth} \label{s:pquple}
\gamma \sim \gamma' \Leftrightarrow (l_{i_1}, \ldots, l_{i_p}, l_{j_1}, \ldots, l_{j_q}) = \lambda (l'_{i_1}, \ldots, l'_{i_p}, l'_{j_1}, \ldots, l'_{j_q}), \quad \text{ o\`u } \lambda \in \C \setminus \{ 0 \}. 
\end{equationth}   
}
\end{definition}

\begin{example}
{\rm 
Consid\'erons l'exemple de la Proposition \ref{profaconetoile}
$$F : \C^2_{(x_1, x_2)} \to \C^2_{(\alpha_1, \alpha_2)}, 
\quad F(x_1,x_2) = \left({(x_1x_2)}^2, {(x_1x_2)}^3 + x_1 \right).$$
Toutes les courbes $\gamma_{\alpha}: (0, + \infty) \to \C^2_{(x_1, x_2)}$ 
tendant vers l'infini telles que $F \circ \gamma_{\alpha}$ 
tend vers un point $(\alpha^2, \alpha^3)$, o\`u $\alpha \neq 0$, sont \'equivalentes. 
Une courbe $\gamma_0$ tendant vers l'infini telle que $F \circ \gamma_0$ tend vers l'origine n'est pas \'equivalente \`a une courbe $\gamma_{\alpha}$, o\`u $\alpha \neq 0$.  
}
\end{example}

Avec ces d\'efinitions, nous pouvons  formaliser maintenant la d\'efinition de fa{\c c}ons \'etoile.

\subsection{Fa{\c c}ons \'etoile}

La Proposition \ref{remarkfaconetoile} ci-dessous est fondamentale. Elle nous permettra de d\'efinir une sous-partition de la partition de $S_F$ d\'efine par les fa{\c c}ons, laquelle sous-partition se r\'ev\'elera \^etre une ``bonne'' stratification. L'id\'ee de la Proposition \ref{remarkfaconetoile}  est la suivante~: \'Etant donn\'ee
  une fa{\c c}on fix\'ee $\kappa$ d'un \'el\'ement $S_\nu$ de la partition de $S_F$ d\'efinie par la relation (\ref{partition facon}), nous subdivisons  $S_\nu$ en utilisant la relation d'\'equivalence  (\ref{s:pquple}) entre  les courbes $\gamma: (0, + \infty) \rightarrow \C^n_{(x)}$ admettant la fa{\c c}on $\kappa$. L'image (par l'application $F$) de 
chaque classe d'\'equivalence   de la relation (\ref{s:pquple})   d\'etermine une sous-vari\'et\'e de $S_\nu$, 
laquelle a une dimension d'autant plus petite que la vitesse de tendre vers z\'ero des courbes coordonn\'ees est grande. Nous obtenons ainsi une 
 sous-partition de  
 $S_{\nu}$. 
\`A partir de ces sous-partitions, nous avons une nouvelle, ``bonne'' partition de $S_F$ (D\'efinition \ref{defXi*}). C'est-\`a-dire, cette fois, cette partition est une stratification (Th\'eor\`eme \ref{theostraXi*}). Le proc\'ed\'e de classer les courbes tendant vers l'infini dans la source, qui sera formalis\'e sous forme de ``fa{\c c}ons \'etoile'' (D\'efinition  \ref{defpreetoile}) est tr\`es significatif~: chaque classe d'\'equivalence de la relation  (\ref{s:pquple}) contient des courbes parall\`eles localement, qui autrement dit, d\'efinissent un feuilletage de dimension (complexe) 1 dans $\C^n_{(x)}$ ({\it cf.}  Exemple \ref{exfeuilletage}). Ce fait est la cl\'e de la d\'emonstration du r\'esultat principal de cet article~: le Th\'eor\`eme \ref{theostraXi*}.

La D\'efinition suivante est important dans toute la suite. 

\begin{definition} \label{generique}
{\rm Soient $S_\nu$ une sous-vari\'et\'e de $S_F$  et $\kappa$ une fa{\c c}on de $S_\nu$. Soit $a$ un point de $S_\nu$ et $\gamma : (0, +\infty) \, \to \C^n_{(x)}$ une courbe tendant vers l'infini avec la fa{\c c}on $\kappa$ telle que $F \circ \gamma$ tende vers $a$. 
Soit $\{ \hat{\gamma}_i \}$ l'ensemble des courbes \'equivalentes \`a $\gamma$ ({\it cf.} D\'efinition \ref{defs:pquple}). 
Nous disons que le point $a$ est {\it un point g\'en\'erique de  $S_\nu$ relativement \`a la fa{\c c}on $\kappa$} s'il 
admet  un voisinage ouvert $U_a$ de $a$ dans  $S_\nu$ et il existe un sous-ensemble des courbes $\{ \hat{\gamma}_j \} \subset \{ \hat{\gamma}_i \}$  telle que l'ensemble des limites de $F \circ \hat{\gamma_j}$  est $U_a$. 

}
\end{definition}

\begin{example}  \label{exemplegenerique}
{\rm Revenons \`a l'exemple $F(x_1,x_2) = \left({(x_1x_2)}^2, {(x_1x_2)}^3 + x_1 \right)$  de la Proposition \ref{profaconetoile}. 
Consid\'erons $S_F$ comme une sous-vari\'et\'e de $S_F$ lui-m\^eme.
 Un point $a = (\alpha^2, \alpha^3)$, o\`u $\alpha \neq 0$, 
est un point g\'en\'erique de $S_F$ relativement \`a la fa{\c c}on $\kappa =(2)[1]$, puisqu'il existe un voisinage $U_a = \{ b = (\beta^2, \beta^3), \beta > 0 \}$ du point $a$ tel que : 

+ toutes les courbes $\gamma_{\beta} = ( 1/u, \beta u)$ 
sont \'equivalentes \`a la courbe  $\gamma_{\alpha} = ( 1/u, \alpha u)$, 

+ l'ensemble des limites $F \circ \gamma_\beta$  est $U_a$. 
}
\end{example}

\begin{remark} \label{remarkgeneric}
{\rm L'ensemble des points g\'en\'eriques de  $S_\nu$ relativement \`a une fa{\c c}on $\kappa$  dans la D\'efinition \ref{generique} est dense dans $S_\nu$.
}
\end{remark}

\begin{proposition} \label{remarkfaconetoile}
{\rm Soit $S_\nu$ un \'el\'ement de dimension $\nu$ de la partition de $S_F$ d\'efinie par la relation (\ref{partition facon}). Pour chaque fa{\c c}on $\kappa$ de $S_\nu$, la relation d'\'equivalence  (\ref{s:pquple}) 
 nous fournit une partition finie ${\{ S^{\kappa}_{\nu_i}}\}_{i = 0, \ldots, t}$ de $S_\nu$,  o\`u $t \leq \nu$, $\, t \in \N$, telle que 
\begin{enumerate}

\item[1)] $\nu = \dim S^{\kappa}_{\nu_0} > \dim S^{\kappa}_{\nu_1} > \dim S^{\kappa}_{\nu_2} > \cdots > \dim S^{\kappa}_{\nu_t},$ 

\item[2)]$S^{\kappa}_{\nu_i} \cap S^{\kappa}_{\nu_j} = \emptyset, \text{ pour }  0 \leq i, j \leq t   \text{ et }  i \ne j, $

\item[3)] $S^{\kappa}_{\nu_i} \subset \overline{S^{\kappa}_{\nu_j}} \text{ pour } i > j \text{ et } 0 \leq i, j \leq t.$
\end{enumerate}

La partition ${\{ S^{\kappa}_{\nu_i}}\}_{i = 0, \ldots, t}$ est appel\'ee
 {\it partition de $S_\nu$ d\'efinie par la fa{\c c}on $\kappa$}.

}
\end{proposition}

\begin{preuve}
Consid\'erons $S_\nu$ un \'el\'ement de dimension $\nu$ de la partition de $S_F$ d\'efinie par la relation (\ref{partition facon}) 
et consid\'erons une fa{\c c}on $\kappa$ de $S_\nu$. D'apr\`es la Remarque \ref{remarkgeneric}, 
l'ensemble des points g\'en\'eriques de $S_\nu$  relativement \`a la fa{\c c}on $\kappa$ ({\it cf.}  D\'efinition  \ref{generique}), not\'e $S^{\kappa}_{\nu_0} $,  
est dense dans $S_\nu$. Nous avons 
 $$ \dim S^{\kappa}_{\nu_0} = \dim S_{\nu}, \quad \text{ et } \quad \overline{S^{\kappa}_{\nu_0}} = S_{\nu}.$$
Si $S^{\kappa}_{\nu_0} = S_{\nu}$, alors la d\'emonstration est finie. Sinon, nous r\'ep\'etons le proc\'ed\'e ci-dessus~: Notons $S^{\kappa}_{\nu_1}$ l'ensemble des points g\'en\'eriques de $S_{\nu} \setminus S^{\kappa}_{\nu_0}$. Nous avons 
 $$ \dim S^{\kappa}_{\nu_1} = \dim (S_{\nu} \setminus S^{\kappa}_{\nu_0}), \quad \text{ et } \quad \overline{S^{\kappa}_{\nu_1}} = S_{\nu} \setminus S^{\kappa}_{\nu_0}.$$
De plus, nous avons 
$$\quad S^{\kappa}_{\nu_0} \cap S^{\kappa}_{\nu_1} = \emptyset \quad \text{ et } \quad S^{\kappa}_{\nu_1} \subset \overline{S^{\kappa}_{\nu_0}}.$$
Puisque $S^{\kappa}_{\nu_0}$ est dense dans $S_{\nu}$, il vient $\dim (S_{\nu} \setminus S^{\kappa}_{\nu_0}) < \dim S^{\kappa}_{\nu_0}$, nous avons alors 
$$\dim S^{\kappa}_{\nu_0} > \dim S^{\kappa}_{\nu_1}.$$
Si $S^{\kappa}_{\nu_1} = S^{\kappa}_{\nu_0}$, la d\'emonstration est finie. Sinon, nous continuons ce proc\'ed\'e. Puisque $\nu$ est fini,   il existe $t \leq \nu$, $\, t \in \N \setminus \{ 0 \}$ tel que  
les points g\'en\'eriques de  $S^{\kappa}_{\nu_t} $ est  $S^{\kappa}_{\nu_t} $. 
 Donc ${\{ S^{\kappa}_{\nu_i}}\}_{i = 0, \ldots ,t}$ d\'ecrit une partition finie de $S_\nu$. 
\end{preuve}

\begin{example}  \label{exemplesubdivisionetoile}
{\rm 
Revenons \`a l'exemple $F(x_1,x_2) = \left({(x_1x_2)}^2, {(x_1x_2)}^3 + x_1 \right)$  de la Proposition \ref{profaconetoile}. D'apr\`es l'exemple  \ref{exemplegenerique}, 
la partition de $S_F$ d\'efinie par la fa{\c c}on (2)[1], d\'etermin\'ee comme dans la Propostion \ref{remarkfaconetoile}  est $S_F \supset \{ 0 \} \supset \emptyset.$
}
\end{example}

\begin{example}  \label{exemplesubdivisionetoile2}
{\rm Soit $F : \C^3_{(x_1, x_2, x_3)} \to \C^3_{(\alpha_1, \alpha_2, \alpha_3)}$ une application polynomiale dominante telle que 
$$F(x_1, x_2, x_3) = \left( (x_1x_2)^2,  (x_2 x_3)^2, x_1x_2^2x_3 + x_2 \right).$$ 
Nous voyons que l'ensemble asymptotique $S_F$ est le c\^one quadratique 
 $\alpha_3^2 = \alpha_1 \alpha_2$. En fait, choisissons la courbe 
$$\gamma = \left( \alpha u, \frac{1}{u}, \beta u \right)$$
tendant vers l'infini, alors $F \circ \gamma (u)$ tend vers un point $(\alpha^2, \beta^2, \alpha  \beta)$. Les courbes $\gamma$ admettent la fa{\c c}on $\kappa = (1, 3)[2]$. Maintenant :

+ Consid\'erons $S_F$ comme une sous-vari\'et\'e de $S_F$ lui-m\^eme. Tout point $a \in S_F \setminus (0 \alpha_1 \cup 0 \alpha_2)$ est un point g\'en\'erique de $S_F$   relativement \`a la fa{\c c}on $\kappa = (1, 3)[2]$ ({\it cf.} D\'efinition \ref{generique}). En fait, tout les courbes $\gamma = \left( \alpha u, \frac{1}{u}, \beta u \right)$ o\`u $\alpha \neq 0$ et $\beta  \neq 0$ sont \'equivalentes et l'ensemble des limites de leurs images 
par $F$ est $S_F \setminus (0 \alpha_1 \cup 0 \alpha_2)$. 

+ Consid\'erons $S_1 =0 \alpha_1 \cup 0 \alpha_2$. Tout point $a \in S_1 \setminus \{ 0 \}$ est un point g\'en\'erique de $S_1$ (relativement \`a la fa{\c c}on $\kappa = (1, 3)[2]$). En fait, pour tout point $a \in 0 \alpha_1 \setminus \{0\}$, il existe des courbes \'equivalentes $\gamma = \left( \alpha u^2, \frac{1}{u^2}, u \right)$, o\`u $\alpha \neq 0$, telles que l'ensemble des limites de leurs images 
par $F$ est $\alpha_1 \setminus \{0\}$. 
De la m\^eme mani\`ere, pour tout point $a \in 0 \alpha_2 \setminus \{0\}$, il existe des courbes \'equivalentes $\gamma = \left( u, \frac{1}{u^2}, \beta u^2 \right)$, o\`u $\beta \neq 0$, telles que 'ensemble des limites de leurs images 
par $F$ est  $0\alpha_2 \setminus \{0\}$.

La partition de $S_F$ relativement \`a la fa{\c c}on $\kappa = (1, 3)[2]$ au sens de la Proposition  \ref{remarkfaconetoile} est :
$$S_F = \{ \alpha_3^2 = \alpha_1 \alpha_2 \} \supset 0\alpha_1 \cup 0\alpha_2 \supset \{ 0 \}.$$

}
\end{example}

\begin{definition} \label{defpreetoile}
{\rm Soit $S_\nu$ un \'el\'ement de la partition de $S_F$ d\'efinie par la relation (\ref{partition facon})  et soit  ${\{ S^{\kappa}_{\nu_i}}\}_{i = 0, ...,t}$  la partition de $S_\nu$ d\'efinie par une fa{\c c}on $\kappa$ de $S_\nu$ comme dans la Proposition \ref{remarkfaconetoile}. Si $t \geq 1$, 
nous d\'efinissons  les ``{\it fa{\c c}ons \'etoile}'' de la fa{\c c}on $\kappa$, comme suit~: 

la fa{\c c}on $\kappa$ de $ S^{\kappa}_{\nu_1}$ appel\'ee la {\it fa{\c c}on \'etoile} $\kappa^{1*}$,

la fa{\c c}on de $\kappa$ de $ S^{\kappa}_{\nu_2}$ appel\'ee la {\it fa{\c c}on \'etoile} $\kappa^{2*}$,

$\quad \quad \quad \quad \quad \quad \quad \quad$ ...

la fa{\c c}on $\kappa$ de $ S^{\kappa}_{\nu_t}$ appel\'ee la {\it fa{\c c}on \'etoile} $\kappa^{t*}$. 

Par convention, nous disons que la fa{\c c}on $\kappa$ est la fa{\c c}ons $\kappa^{0*}$ de $S_\nu$.
}
\end{definition}

\begin{example} \label{exemplesubdivisionetoile1}
{\rm D'apr\`es l'exemple \ref{exemplesubdivisionetoile}, la fa{\c c}on (2)[1] de l'ensemble asymptotique de l'application $F(x_1,x_2) = \left({(x_1x_2)}^2, {(x_1x_2)}^3 + x_1 \right)$  
  admet une seule fa{\c c}on \'etoile ${(2)[1]}^{1*}$.
}
\end{example}

\begin{example} 
{\rm La fa{\c c}on $\kappa = (1, 3)[2]$ de l'application $F : \C^3_{(x_1, x_2, x_3)} \to \C^3_{(\alpha_1, \alpha_2, \alpha_3)}$ admet deux fa{\c c}ons \'etoile: 

+ fa{\c c}on \'etoile $(1, 3)[2]^{1*}$ correspondante \`a la strate $0 \alpha_1 \cup 0 \alpha_2 \setminus \{ 0\}$,

+ fa{\c c}on \'etoile $(1, 3)[2]^{2*}$ correspondante  \`a l'origine
\\ ({\it cf.} Exemple  \ref{exemplesubdivisionetoile2}).
}
\end{example}

\begin{remark}
{\rm 
Chaque \'el\'ement  $S^{\kappa}_{\nu_i}$ de la partition ${\{ S^{\kappa}_{\nu_i}}\}_{i = 0, \ldots, t}$  d\'efini dans la Proposition \ref{remarkfaconetoile} peut associer aux plusieurs classes d'\'equivanlence des courbes de la relation (\ref{s:pquple}). Par exemple, dans l'exemple 
 \ref{exemplesubdivisionetoile2}, pour tout point $a \in 0 \alpha_1 \setminus \{0\}$, il existe les courbes \'equivalentes $\gamma = \left( \alpha u^2, \frac{1}{u^2}, u \right)$, o\`u $\alpha \neq 0$, mais il existe aussi des autres classes des coubres \'equivalentes, par exemple, $\hat{\gamma} = \left( \alpha u^3, \frac{1}{u^3}, u \right)$, o\`u $\alpha \neq 0$, tendant vers l'infini telles que $F \circ \hat{\gamma}$ tendent vers $a$.
 Alors, dans la suite de cet article, nous utilisons la  convention suivante :  dans un voisinage ouvert  d'un \'el\'ement  $S^{\kappa}_{\nu_i}$ de la partition ${\{ S^{\kappa}_{\nu_i}}\}_{i = 0, \ldots, t}$, nous fixons seulement {\it une} classe d'\'equivalence des courbes correspondantes. 
 Avec cette convention, nous avons :
}
\end{remark}

\begin{lemma} \label{lemmeequivalente}
{\rm Toutes les courbes correspondantes aux points d'un ensemble  ouvert  de chaque \'el\'ement $S^{\kappa}_{\nu_i}$ de la partition ${\{ S^{\kappa}_{\nu_i}}\}_{i = 0, \ldots, t}$ d\'efinie dans la Proposition \ref{remarkfaconetoile} sont \'equi\-va\-lentes.
}
\end{lemma}

Avec le Lemme  \ref{lemmeequivalente}, nous avons la Propostion suivante, qui montre  la signification g\'eom\'etrique des  fa{\c c}ons \'etoile.

\begin{proposition} \label{geometrie faconetoile}
{\rm 
Fixons une fa{\c c}on \'etoile $\kappa^{i*}$ de $S_F$. Soit $U_a$ un voisinage  ouvert su\-ffisamment petit d'un point $a \in S_F$ tel que tout point de $U_a$ admet la fa{\c c}on \'etoile $\kappa^{i*}$. 
Alors toutes les courbes correspondantes aux points de $U_a$  d\'efinissent un feuilletage de dimension (complexe) 1 de $\C^n_{(x)}$.

}
\end{proposition}

\begin{figure}[h!]
\centering
\includegraphics[scale=0.7]{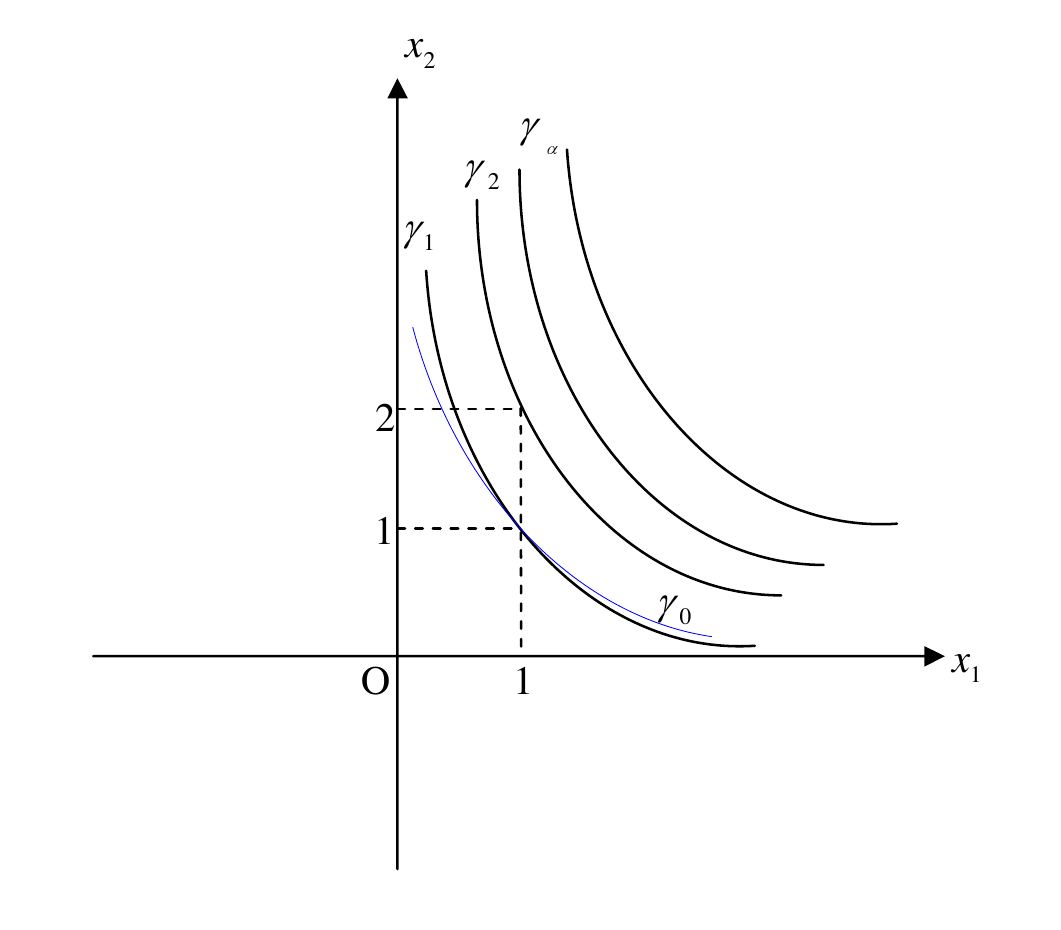}
\caption{Feuilletage des courbes de la m\^eme classe d'\'equivalence (\ref{s:pquple}).}
\label{facon2}
\end{figure}

\begin{preuve} 
{\rm Soit $U_a$  un voisinage  ouvert suffisamment petit d'un point $a \in S_F$ 
 tel que tout point de $U_a$ admet la fa{\c c}on \'etoile $\kappa^{i*}$. 
 Soient $b, b' \in U_a$  et  $\gamma, \gamma':  (0, +\infty) \, \to \C^n_{(x)}$  deux courbes tendant vers l'infini  telles que $F \circ \gamma$ et  $F  \circ \gamma'$ tendent vers 
$b$ et $b'$, respectivement. D'apr\`es  la  D\'efinition \ref{defpreetoile},  comme $b$, $b'$ admettent la m\^eme fa\c con \'etoile $\kappa^{i*}$,  alors ils appartiennent au m\^eme \'el\'ement $S_{\nu_i}^\kappa$ d\'efini dans la  Proposition \ref{remarkfaconetoile}. 
D'apr\`es  le Lemme \ref{lemmeequivalente}, deux courbes $\gamma$ et $\gamma'$ sont \'equivalentes. D'apr\`es la D\'efinition  \ref{defs:pquple}, 
les $(p+q)$-uples  
associ\'es aux deux courbes correspondantes $\gamma$ et $\gamma'$
 sont proportionnels. 
Donc les deux vecteurs tangents $d \gamma (u)$, $d\gamma' (u)$ \`a  ces deux courbes en une m\^eme valeur de param\`etre $u \in (0, +\infty) $ sont parall\`eles.   L'ensemble des courbes correspondantes aux points de $S_\nu$ admettant la m\^eme fa\c con \'etoile $\kappa^{i*}$ d\'efinit un feuilletage de dimension (complexe) 1 de $\C^n_{(x)}$ (voir Figure \ref{facon2}).

}
\end{preuve}

\begin{example}  \label{exfeuilletage}
{\rm Prenons l'exemple 
$$F : \C^2_{(x_1, x_2)} \to \C^2_{(\alpha_1, \alpha_2)}, \quad F(x_1,x_2) = \left({(x_1x_2)}^2, {(x_1x_2)}^3 + x_1 \right)$$
({\it cf.} Proposition \ref{profaconetoile}) pour illustrer la Proposition \ref{geometrie faconetoile}. Nous savons que l'ensemble asymptotique $S_F$ de l'application $F$ est le {\it cusp} $\{(\alpha^2, \alpha^3)~: \alpha \in \C \}$. Prenons $\alpha \neq 0$, par exemple~:

- Si $\alpha = 1$, alors il existe la courbe 
$$\gamma_1 : (0, +\infty) \, \to \C^2_{(x_1, x_2)}, \quad \gamma(u) = \left( \frac{1}{u}, u \right)$$
tendant vers l'infini telle que $F \circ \gamma_1$ tend vers le point $(1, 1) \in S_F$.

- Si $\alpha = 2$, alors il existe la courbe ~: 
$$\gamma_2 : (0, +\infty) \, \to \C^2_{(x_1, x_2)}, \quad \gamma(u) = \left( \frac{1}{u}, 2u \right)$$
tendant vers l'infini telle que $F \circ \gamma_2$ tend vers le point $(4, 8) \in S_F$. 

Nous voyons que deux courbes  $\gamma_1$ et $\gamma_2$ sont parall\`eles localement. 
 En g\'en\'eral, toutes les courbes $\gamma_\alpha$ tendant vers l'infini telles que $F \circ \gamma_\alpha$ tendent vers $(\alpha^2, \alpha^3) \in S_F$ o\`u $\alpha \neq 0$  forment un ensemble ${ \{\gamma_\alpha \}}_{\alpha \neq 0}$ de courbes parall\`eles localement. C'est-\`a-dire, ${ \{\gamma_\alpha \}}_{\alpha \neq 0}$ est un feuilletage de $\C^2_{(x_1, x_2)}$ (voir figure \ref{facon2}). 

Par contre, quelque soit une courbe correspondant \`a l'origine, par exemple, $\gamma_0 = \left\{ \frac{1}{u^2}, u \right\}$ n'est pas ``parall\`ele'' avec la courbe $\gamma_\alpha$, pour tout $\alpha \neq 0$ (voir figure \ref{facon2}).
 }
\end{example}

\begin{remark}
{\rm D'apr\`es la Lemme \ref{lemmeequivalente}, nous voyons que pour chaque  fa{\c c}on \'etoile $\kappa^{i*}$ de $S_F$, $i \geq 0$, alors avec un changement de param\`etre $u$, il existe une unique courbe $\gamma(u) : (0, +\infty) \, \to \C^n_{(x)}$ tendant vers l'infini  avec la fa{\c c}on \'etoile $\kappa^{i*}$ et telle que $F \circ \gamma (u)$ tend vers $a \in S_F$.
 Par exemple, soit $F : \C^2_{(x_1, x_2)} \to \C^2_{(\alpha_1, \alpha_2)}$ 
l'application polynomiale telle que   $F(x_1,x_2) = (x_1, x_1x_2)$. 
 L'ensemble asymptotique $S_F$ est le plan $\alpha_1 = 0$ dans $\C^2_{(\alpha_1, \alpha_2)}$. Prenons le point $(0,1) \in S_F$. Il existe des suites 
 $\left\{ \left( \frac{1}{k}, k \right) \right\}, 
\left\{ \left( \frac{1}{k^2}, k^2 \right) \right\}, 
\left\{ \left( \frac{1}{k^3}, k^3 \right) \right\} \ldots $ 
tendant vers l'infini telles que leurs images tendent vers le point $(0, 1) \in S_F$. En fait, l'image de toute suite  $\left\{ \left( \frac{1}{\alpha k^r}, \alpha k^r \right) \right\}$, o\`u $r \in \N \setminus \{0\}$ et $\alpha \in \C \setminus \{ 0 \}$, tend vers le point $(0, 1)$. Par un changement de param\`etre, ces suites repr\'esentent la seule courbe  
$$\gamma : (0, +\infty) \, \to \C^2_{(x_1, x_2)}, \quad \gamma(u) = \left( \frac{1}{u}, u \right).$$
}
\end{remark}

\section{Partition de l'ensemble asymptotique $S_F$ d\'efinie par les fa{\c c}ons \'etoile}
Notons $ \Xi^*(a)$ l'ensemble des fa\c cons \'etoile du point $a$. 
Suivant la Proposition \ref{pro partitionfini1}, l'id\'ee naturelle pour stratifier l'ensemble asymptotique $S_F$ est de consid\'erer la 
partition de $S_F$ d\'efinie par la relation: $ a_1 \sim a_2 \text{ si et seulement si } \Xi^*(a_1) = \Xi^*(a_2).$

 Cependant, cette partition contient des ``fausse strates'', elles-m\^emes plong\'ees dans une vari\'et\'e lisse 
de plus grande dimension. Par exemple, un point $a$ de $S_\nu$ peut \^etre un point g\'en\'erique (au sens de la D\'efinition \ref{generique})
 pour 
une fa{\c c}on $\kappa_1$ et n'\^etre pas g\'en\'erique pour une fa{\c c}on $\kappa_2 \neq \kappa_1$. 
Afin d'\'eviter cet inconv\'enient, nous utiliserons la partition d\'efinie comme suit~: 

\begin{definition} \label{defXi*} 
{\rm Consid\'erons  $S_\nu$ un \'el\'ement de la stratification de $S_F$ d\'efinie par les fa{\c c}ons
et $\Xi(S_\nu)$ l'ensemble  de toutes les fa{\c c}ons de $S_\nu$. 
Nous d\'efinissons $S_{\nu_0}$ comme l'ensemble des points g\'en\'eriques relativement \`a  
au moins une fa{\c c}on de $S_\nu$ (au sens de la D\'efinition \ref{generique}), c'est-\`a-dire ~: 
$$S_{\nu_0} = \bigcup_{\kappa \, \in \, \Xi(S_\nu)} S^{\kappa}_{\nu_0},$$
 o\`u $S_{\nu_0}^\kappa$ est d\'efini comme dans la Proposition \ref{remarkfaconetoile}.

 Notons $A_{\nu_0} = S_\nu \setminus S_{\nu_0}$. 
Nous d\'efinissons $S_{\nu_1}$ comme l'ensemble des points g\'en\'eriques  relativement \`a au moins une fa{\c c}on de $A_{\nu_0}$ 
et nous notons $A_{\nu_1} = A_{\nu_0} \setminus S_{\nu_1}$. En g\'en\'eral, pour $i \geq 1$, 
nous d\'efinissons 
$S_{\nu_i}$ comme l'ensemble des points g\'en\'eriques  relativement \`a au moins une fa{\c c}on de $A_{\nu_{i-1}} $ et 
$A_{\nu_i} = A_{\nu_{i-1}} \setminus S_{\nu_i}$. Nous obtenons une   sous-partition de $S_\nu$. 
 \`A partir de ces sous-partitions, nous obtenons ainsi une nouvelle partition de $S_F$, 
appel\'ee {\it la partition de $S_F$ d\'efinie par les fa{\c c}ons \'etoile}. 
}
\end{definition}

\begin{example} 
{\rm Revenons \`a l'exemple $F(x_1,x_2) = \left({(x_1x_2)}^2, {(x_1x_2)}^3 + x_1 \right)$  de la Proposition \ref{profaconetoile}. D'apr\`es l'exemple \ref{exemplesubdivisionetoile1}, nous avons la partition de $S_F$ d\'efinie par les fa{\c c}ons \'etoile est d\'efinie par la filtration $S_F \supset \{ 0 \} \supset \emptyset.$

}
\end{example}

\begin{example} 
{\rm Revenos \`a l'exemple \ref{exfacon}, nous avons:

- La fa{\c c}on $(3)[1]$ n'a pas de fa{\c c}ons \'etoile. 

- La fa{\c c}on $(3)[2]$ n'a pas de fa{\c c}ons \'etoile. 

- La fa{\c c}on $(3)[1, 2]$ n'a qu'une seule fa{\c c}on \'etoile $(3)[1, 2]^{1*}$.

La partition de l'ensemble asymptotique $S_F$ d\'efinie par les fa{\c c}ons \'etoile de l'exemple \ref{exfacon} est 
$$S_F  \supset 0 \alpha_3 \supset \{ 0 \}.$$ 

}
\end{example}

\begin{example} 
{\rm La partition de l'ensemble asymptotique $S_F$ d\'efinie par les fa{\c c}ons \'etoile de l'exemple \ref{exemplesubdivisionetoile2}  est 
$$S_F  \supset 0\alpha_1 \cup 0\alpha_2 \supset \{ 0 \}.$$

}
\end{example} 

\begin{example}  
{\rm Soit l'appication polynomiale dominante
$$F : \C^3_{(x_1, x_2, x_3)} \to \C^3_{(\alpha_1, \alpha_2, \alpha_3)}, \quad F = (x_1 x_2, (x_2 x_3)^2, x_1x_2^2x_3 + x_2).$$
 L'ensemble asymptotics $S_F$ est le parapluie de Whitney $(\alpha, \beta^2, \alpha \beta)$. 
 En fait, choisissons la courbe 
$$\left( \alpha u, \frac{1}{u}, \beta u \right),$$
 de la fa{\c c}on $\kappa = (1, 3) [ 2 ]$, 
tendant vers l'infini, alors $F \circ \gamma (u)$ tend vers un point $(\alpha, \beta^2, \alpha \beta)$. De la m\^eme mani\`ere que dans l'exemple \ref{exemplesubdivisionetoile2}, nous avons 
la partition de $S_F$ relativement \`a la fa{\c c}on $\kappa = (1, 3)[2]$ au sens de la Proposition  \ref{remarkfaconetoile} est 
$$S_F \supset 0\alpha_1 \cup 0\alpha_2 \supset \{ 0 \},$$ et ce-ci est aussi la partition de l'ensemble asymptotique $S_F$ d\'efinie par les fa{\c c}ons \'etoile.

}
\end{example}

\section{Theor\`eme sur une stratification de l'ensemble asymptotique}

\begin{theorem} \label{theostraXi*} 
Soit $F : \C_{(x)}^n \to \C_{(\alpha)}^n$  une application polynomiale dominante. La partition de l'ensemble asymptotique de $F$ d\'efinie par les fa{\c c}ons \'etoile est une stratification  diff\'erentiable satisfaisant la propri\'et\'e de fronti\`ere.
\end{theorem}

\begin{preuve}
Supposons que $({\mathscr{S}})$ soit la partition de $S_F$ d\'efinie par les fa{\c c}ons \'etoile ({\it cf.} D\'efinition \ref{defXi*}).  
  Nous prouvons maintenant que la partition $({\mathscr{S}})$ est une stratification diff\'erentiable. 
Consid\'erons $S_{\nu_i}$ 
une sous vari\'et\'e de dimension $\nu_i$  de $({\mathscr{S}})$. 
Montrons qu'avec  tout point $a$ de  
$S_{\nu_i}$ admet 
un voisinage dans  $S_{\nu_i}$ diff\'eomorphe \`a une boule de dimension  $\nu_i$.  Par la D\'efinition \ref{defXi*}, nous avons 
$$S_{\nu_i} = \bigcup_{\kappa \, \in \, \Xi(S_{\nu_i})} S^{\kappa}_{\nu_i},$$
 o\`u $\Xi(S_{\nu_i})$ est l'ensemble des fa{\c c}ons de $S_{\nu_i}$ e $S^{\kappa}_{\nu_i}$ est d\'efini comme dans la Proposition \ref{remarkfaconetoile}.
Comme $a \in S_{\nu_i}$, alors il existe une fa{\c c}on $\kappa$ de $\Xi(S_{\nu_i})$ telle que $a \in S^{\kappa}_{\nu_i}$. 
Fixons cette fa{\c c}on $\kappa$. 
Comme $S_{\nu_i}^{\kappa}$ est ouvert, il existe un voisinage ouvert $U_a$ de $a$ dans $S_{\nu_i}^{\kappa}$. 
 Pour tout point $a' \in U_a$, il existe une courbe  diff\'erentible $\gamma^{a'}(u): (0, +\infty) \to  \C_{(x)}^n$   
tendant vers l'infini et telle que $F \circ \gamma^{a'}(u)$ tend vers $a'$ quand $u$ tend vers l'infini. 
 D'apr\`es la D\'efinition \ref{defpreetoile}, tout point de $U_a$ admet la m\^eme fa{\c c}on $\kappa^{i*}$. 
D'apr\`es la Proposition \ref{geometrie faconetoile}, 
 l'ensemble des courbes $\gamma^{a'}(u)$  telles que $F(\gamma^{a'}(u))$ tende vers $a'$ pour $a' \in U_a$ 
d\'efinit un feuilletage $\mathcal{F}$ de dimension complexe 1 de $\C^n_{(x)}$ et un feuilletage transverse $\mathcal{F'}$ de dimension complexe $(n-1)$.  Ce feuilletage transverse est  constitu\'e de sous-vari\'et\'es  transverses aux feuilles locales de $\mathcal{F}$. 
Les espaces tangents aux feuilletages $\mathcal{F}$ et $\mathcal{F'}$ sont engendr\'es res\-pectivement par des champs de vecteurs $v_1^{a'}(u)$ et $v_2^{a'}(u), \ldots, v_{n}^{a'}(u)$ d\'ependant du point $a'$ et du param\`etre $u$. 
Pour chaque valeur de param\`etre $u$, les images des feuilletages $\mathcal{F}$ et $\mathcal{F'}$  sont engendr\'ees respectivement par des champs de vecteurs $dF(v_1^{a'}(u))$ et  $dF(v_2^{a'}(u)),$ $\ldots,$ $dF(v_{n}^{a'}(u))$. 
D'une part, lorsque $u$ tend vers l'infini, $dF(v_1^{a'}(u))$, qui est tangent \`a la courbe $F(\gamma^{a'}(u))$, tend vers 0 puisque la courbe $F(\gamma^{a'}(u))$ tend vers le point $a'$. 
D'autre part, la  signification g\'eom\'etrique de la Proposition \ref{remarkfaconetoile} montre
que plus la dimension de $S_{\nu_i}$ est petite, plus il y a de courbes coordonn\'ees 
de la courbe $\gamma^{a'}(u)$ tendant vers z\'ero plus rapidement et donc 
$(n - 1) - \nu_i$ 
vecteurs parmi les vecteurs $dF(v_2^{a'}(u)),$ $\ldots,$ $dF(v_{n}^{a'}(u))$ 
tendant vers z\'ero quand $u$ tend vers l'infini. La dimension de l'espace tangent \`a $U_a$ au point $a'$ est $\nu_i$. Donc $U_a$ est diff\'eomorphe \`a une boule de dimension $\nu_i$. La partition $({\mathscr{S}})$ est donc une stratification diff\'erentiable. 

La  propri\'et\'e de fronti\`ere de la stratification $({\mathscr{S}})$ d\'ecoule directement de la D\'efinition \ref{defXi*}.
\end{preuve}

\bibliographystyle{plain}

\end{document}